\documentclass[10pt,a4paper]{article}

\usepackage{amssymb}
\usepackage{amsthm}
\usepackage{amsmath}

\newtheorem{thm}{Theorem}

\newtheorem{lem}{Lemma}

\newcommand*{\half}{\frac12}
\newcommand*{\quart}{\frac14}

\newcommand*{\lb}{\left\{}
\newcommand*{\rb}{\right\}}

\newcommand*{\PP}{\mathbf{P}}
\newcommand*{\EE}{\mathbf{E}}
\newcommand*{\Var}{\mathbf{Var}}

\newcommand*{\loc}{\ensuremath{\mathcal{L}}}

\begin{document}

\title{An elementary approach to Brownian local time based on
simple, symmetric random walks}

\author{Tam\'as Szabados\footnote{Corresponding author, address:
Department of Mathematics, Budapest University of Technology and
Economics, M\H{u}egyetem rkp. 3, H \'ep. V em. Budapest, 1521,
Hungary, e-mail: szabados@math.bme.hu, telephone: (+36 1)
463-1111/ext. 5907, fax: (+36 1) 463-1677} \footnote{Research
supported by a Hungarian National Research Foundation (OTKA) grant
No. T42496.} and Bal\'azs Sz\'ekely\footnote{Research
supported by the HSN laboratory of BUTE.} \\
Budapest University of Technology and Economics}

\date{}

\maketitle

\begin{center}
\emph{Dedicated to the 70th birthdays of Endre Cs\'aki and P\'al
R\'ev\'esz}
\end{center}

\bigskip


\begin{abstract}
In this paper we define Brownian local time as the almost sure
limit of the local times of a nested sequence of simple, symmetric
random walks. The limit is jointly continuous in $(t,x)$. The rate
of convergence is $n^{\frac14} (\log n)^{\frac34}$ that is close
to the best possible. The tools we apply are almost exclusively
from elementary probability theory.
\end{abstract}


\renewcommand{\thefootnote}{\alph{footnote}}
\footnotetext{ 2000 \emph{MSC.} Primary 60J55, 60J65. Secondary
60F15.} 

\section{Introduction}

The present work is part of a bigger project that aims to rebuild
stochastic calculus using almost sure (strong) approximations by
simple, symmetric random walks (RW's), applying tools almost
exclusively from elementary probability theory. The most advanced
tool one needs for this project is a large deviation inequality.
The underlying motivation is partly didactic, partly the belief
that this elementary approach may help in attacking some harder
problems. The prototype of such efforts was the construction of
Brownian motion (BM = Wiener process) as an almost sure limit of
simple RW paths, given by Frank Knight in 1962 \cite{Kni62}. This
lead him to a related construction of Brownian local time via the
now celebrated Ray-Knight theory in 1963 \cite{Kni63}. An
important contribution to the theory of local time was given by
P\'al R\'ev\'esz in 1981 \cite{Rev81} in which he showed that
Brownian local time can be almost surely approximated by the local
time of a simple, symmetric RW with rate of convergence
$n^{\frac{1}{4}+\epsilon}$. Further contributions and
generalizations to this were given by Endre Cs\'aki and P\'al
R\'ev\'esz in 1983 \cite{CsR83}, A.N. Borodin in 1989
\cite{Bor89}, and several others. M. Cs\"org\H{o} and L. Horv\'ath
in 1989 \cite{CsH89} gave a best possible strong approximation
based on Skorohod embedding with rate $n^{\frac14} (\log
n)^{\frac12} (\log \log n)^{\frac14}$. This was generalized for a
wide class of RW's by R.F. Bass and D. Khoshnevisan in 1993
\cite{BKh93}.

The level sets $\{t: W(t)=x, \: 0 \le t < \infty \}$ of Brownian
motion have zero Lebesgue measure, hence it is not at all obvious
that Brownian local time given by the classical definition of P.
L\'evy
\begin{equation} \label{eq:Levyloctime}
\loc(t,x) = \lim_{\epsilon \searrow 0} \frac{1}{2 \epsilon} \:
\lambda \{s \in [0,t] : W(s) \in (x-\epsilon , x+\epsilon) \} ,
\end{equation}
where $\lambda$ denotes Lebesgue measure, is a well-defined,
non-vanishing process which has a version jointly continuous in
$(t,x)$. This fact was first shown by Trotter in 1958
\cite{Tro58}. In this paper we define Brownian local time as the
almost sure limit of the local times of a nested sequence of
simple, symmetric RW's. This limit is automatically jointly
continuous in $(t,x)$. The rate of convergence is $n^{\frac14}
(\log n)^{\frac34}$ that is close to the best possible.

\section{From random walks to Brownian motion}

A main tool of the present paper is an elementary construction of
BM. The specific construction we are going to use in the sequel,
taken from \cite{Sza96}, is based on a nested sequence of simple,
symmetric random walks that uniformly converges to the Wiener
process on bounded intervals with probability $1$. This will be
called \emph{``twist and shrink''} construction or RW construction
in the sequel.

We summarize the major steps of the \emph{``twist and shrink''}
construction here. We start with an infinite matrix of independent
and identically distributed (i.i.d.) random variables $X_m(k)$,
$\PP \lb X_m(k)= \pm 1 \rb = 1/2$ ($m\ge 0$, $k\ge 1$), defined on
the same underlying probability space $(\Omega,\mathcal{F},\PP)$.
Each row of this matrix is a basis of an approximation of the
Wiener process with a dyadic step size $\Delta t=2^{-2m}$ in time
and a corresponding step size $\Delta x=2^{-m}$ in space: $S_m(0)
= 0$, $S_m(n) = \sum_{k=1}^{n} X_m(k)$ $(n \ge 1)$.

The second step of the construction is \emph{twisting}. From the
independent RW's we want to create dependent ones so that after
shrinking temporal and spatial step sizes, each consecutive RW
becomes a refinement of the previous one.  Since the spatial unit
will be halved at each consecutive row, we define stopping times
by $T_m(0)=0$, and for $k\ge 0$,
\[
T_m(k+1)=\min \{n: n>T_m(k), |S_m(n)-S_m(T_m(k))|=2\} \qquad (m\ge
1)
\]
These are the random time instants when a RW visits even integers,
different from the previous one. After shrinking the spatial unit
by half, a suitable modification of this RW will visit the same
integers in the same order as the previous RW. We operate here on
each point $\omega\in\Omega$ of the sample space separately, i.e.
we fix a sample path of each RW. We define twisted RW{}s
$\tilde{S}_m$ recursively for $k=1,2,\dots$ using
$\tilde{S}_{m-1}$, starting with $\tilde{S}_0(n)=S_0(n)$ $(n\ge
0)$ and $S_m(0) = 0$ for any $m \ge 0$. With each fixed $m$ we
proceed for $k=0,1,2,\dots$ successively, and for every $n$ in the
corresponding bridge, $T_m(k)<n\le T_m(k+1)$. Any bridge is
flipped if its sign differs from the desired:
\[
\tilde{X}_m(n)=\left\{
\begin{array}{rl}
 X_m(n)& \mbox{ if } S_m(T_m(k+1)) - S_m(T_m(k))
= 2\tilde X_{m-1}(k+1), \\
- X_m(n)& \mbox{ otherwise,}
\end{array}
\right.
\]
and then $\tilde{S}_m(n)=\tilde{S}_m(n-1)+\tilde{X}_m(n)$. Then
$\tilde{S}_m(n)$ $(n\ge 0)$ is still a simple symmetric random
walk \cite[Lemma 1]{Sza96}. The twisted RW's have the desired
refinement property:
\[
\tilde{S}_{m+1}(T_{m+1}(k)) = 2 \tilde{S}_{m}(k) \qquad (m\ge 0,
k\ge 0).
\]

The third step of the RW construction is \emph{shrinking}. The
sample paths of $\tilde{S}_m(n)$ $(n\ge 0)$ can be extended to
continuous functions by linear interpolation, this way one gets
$\tilde{S}_m(t)$ $(t\ge 0)$ for real $t$. Then we define the $mth$
\emph{approximating RW} by
\[
\tilde{B}_m(t)=2^{-m}\tilde{S}_m(t2^{2m}).
\]
Then the \emph{refinement property} takes the form
\begin{equation}
\tilde{B}_{m+1}\left(T_{m+1}(k)2^{-2(m+1)}\right) = \tilde{B}_m
\left( k2^{-2m}\right) \qquad (m\ge 0,k\ge 0). \label{eq:refin}
\end{equation}
Note that a refinement takes the same dyadic values in the same
order as the previous shrunken walk, but there is a \emph{time
lag} in general:
\begin{equation} T_{m+1}(k)2^{-2(m+1)} - k2^{-2m} \ne 0 .
\label{eq:tlag}
\end{equation}

Then we quote some important facts from \cite{Sza96} about the
above RW construction that will be used in the sequel. These will
be stated in somewhat stronger forms but can be read easily from
the proofs in the cited reference, cf. Lemmas 2-4 and Theorem 3
there.

\begin{lem} \label{le:ldi}
Suppose that $X_1,X_2,\dots $ is a sequence of i.i.d. random
variables, $\EE(X_k)=0$, $\Var(X_k)=1$, and their moment
generating function is finite in a neighborhood of 0. Let $S_k =
X_1 + \dots + X_k$, $k \ge 1$. Then for any $C>1$ and $N \ge
N_0(C)$ one has
\[
\PP \lb \sup_{1 \le k \le N} |S_k| \ge (2CN \log N)^{\half} \rb
\le 2N^{1-C} .
\]
\end{lem}
We mention that this basic fact, that appears in the
above-mentioned reference \cite{Sza96}, essentially depends on a
large deviation theorem. Also, the lemma is valid even when $N$ is
not an integer.

We have a more convenient result in a special case of Hoeffding's
inequality. Let $S_1, S_2, \dots , S_N$ be arbitrary (finite or
countable) sums of the type: $S_k = \sum_r a_{k r} X_{k r}$ with
$\PP \lb X_{k r} = \pm 1 \rb = \half$, $1 \le k \le N$, where
$X_{k r}$ and $X_{l s}$ can be dependent when $k \ne l$. Then we
have the following analog of Lemma~\ref{le:ldi}, cf. \cite{Sza01}
or \cite{SzeSza04}: for any $C >1$ and $N \ge 1$,
\begin{equation}
\PP \lb \sup_{1\le k \le N} |S_k| \ge \left(2C\log N \: \sup_{1\le
k \le N} \Var (S_k)\right) ^\half \rb \le 2 N^{1-C}.
\label{eq:ldp2}
\end{equation}

Lemma \ref{le:ldi} easily implies that the time lags
(\ref{eq:tlag}) are uniformly small if $m$ or $K$ are large
enough.
\begin{lem} \label{le:tlag}
For any $C>1$, and for any $K>0$ and $m \ge 1$ such that $K 2^{2m}
\ge N_1(C)$, we have
\begin{eqnarray*}
\PP \lb \sup_{0\le k2^{-2m}\le K} |T_{m+1}(k)2^{-2(m+1)}
-k2^{-2m}| \ge \left(\frac32 C K \log{}_*K\right)^{\half}
m^{\half} 2^{-m} \rb \\
\le 2(K2^{2m})^{1-C} ,
\end{eqnarray*}
where $\log_* x = \max \{ 1,\log x \}$.
\end{lem}

This lemma and the refinement property (\ref{eq:refin}) implies
the uniform closeness of two approximations of BM if $m$ or $K$
are large enough. We give a proof of the next lemma since its
statement is somewhat stronger than the one of the corresponding
Lemma 4 in \cite{Sza96}; also, we want to emphasize the
similarities between this and lemmas about the local time below.

\begin{lem} \label{le:refin}
For any $C>1$, and for any $K>0$ and $m \ge 1$ such that $K 2^{2m}
\ge N_2(C)$, we have

(a)
\begin{eqnarray*}
\PP \lb \sup_{0\le k2^{-2m} \le K} |\tilde{B}_{m+1}(k2^{-2m})
-\tilde{B}_m(k2^{-2m})| \ge \frac{11}{4} C K^{\quart}
(\log_*K)^{\frac34}
m^{\frac34}  2^{-\frac{m}{2}} \rb \\
\le 6 (K2^{2m})^{1-C} ,
\end{eqnarray*}

(b)
\begin{eqnarray*}
\PP \lb \sup_{j \ge 1} \sup_{0 \le t \le K} |\tilde{B}_{m+j}(t)
-\tilde{B}_m(t)| \ge 27 \: C K_*^{\quart} (\log_*K)^{\frac34}
m^{\frac34}  2^{-\frac{m}{2}} \rb \\
\le \frac{6}{1-4^{1-C}} (K2^{2m})^{1-C} ,
\end{eqnarray*}

where $K_* = \max\{1, K\}$.
\end{lem}

\begin{proof}

(a) Using the abbreviation $t_k = k2^{-2m}$, we have
\begin{eqnarray*}
\tilde B_{m+1}(t_k) -  \tilde B_m(t_k) & = & \tilde B_{m+1}\left(
4k 2^{-2(m+1)}\right) -
\tilde B_{m+1}\left( T_{m+1}(k) 2^{-2(m+1)}\right)  \\
& = & 2^{-m-1}\left(\tilde S_{m+1}(4k) - \tilde
S_{m+1}(T_{m+1}(k))\right) .
\end{eqnarray*}
Let $D_{K,m} = C K^{\quart} (\log_*K)^{\frac34} m^{\frac34}
2^{-\frac{m}{2}}$. Then
\begin{eqnarray*}
\lefteqn{\PP \lb \sup_{0\le t_k \le K} |\tilde{B}_{m+1}(t_k)
-\tilde{B}_m(t_k)| \ge \frac{11}{4} D_{K,m} \rb} \\
& \le & \PP \lb A_{K,m} \rb + \sum\limits_{k=1}^{K 2^{2m}} \PP \lb
\sup_{\{j: |j-4k| \le N'\}} \left| \tilde S_{m+1}(j)
 - \tilde S_{m+1}(4k) \right|
\ge \frac{11}{2} 2^m D_{K,m} \rb  ,
\end{eqnarray*}
where
\begin{equation}\label{eq:AKm}
A_{K,m} = \lb \sup_{0 \le k 2^{-2m} \le K } |T_{m+1}(k) - 4k| \ge
N' \rb , \quad N' = \lfloor(24 C m K \log_* K)^\half 2^m\rfloor .
\end{equation}

Here we can apply Lemma \ref{le:tlag} to the first term and
inequality (\ref{eq:ldp2}) with $N'$ and $C'$ to the terms in the
summation above. The constant $C'>1$ is chosen so that the two
terms be of the same order: $(N')^{1-C'} \approx N^{-C}$, where
$N=K 2^{2m}$. Then $(N')^{1-C'} \le (K^{\half} 2^m)^{1-C'} = (K
2^{2m})^{-C}$ implies $C' = 2C + 1 < 3C$. Since $\log N' \le m
\log_*C \log_*K$ if $K$ or $m$ are large enough and
\[
\sup_{\{j: |j-4k| \le N'\}} \Var\left(\tilde S_{m+1}(j) - \tilde
S_{m+1}(4k)\right) = N' ,
\]
it follows that $(2C' N' \log N')^{\half} \le \frac{11}{2} 2^m
D_{K,m}$.

As a result, we obtain that
\begin{eqnarray*}
& \PP & \lb \sup_{0\le t_k \le K} |\tilde{B}_{m+1}(t_k)
-\tilde{B}_m(t_k)| \ge \frac{11}{4} D_{K,m} \rb \\
&& \le 2 (K 2^{2m})^{1-C} + 2 K 2^{2m} 2 (N')^{1-C'} \le 6 (K
2^{2m})^{1-C} .
\end{eqnarray*}
This proves (a).

(b) Let $D^*_{K,m} = C K_*^{\quart} (\log_*K)^{\frac34}
m^{\frac34} 2^{-\frac{m}{2}}$. By (a),
\[
\sup_{0\le t_k \le K} |\tilde B_{m+1}(t_k) -  \tilde B_m(t_k)| <
\frac{11}{4} D^*_{K,m} ,
\]
except for an event of probability not exceeding $6 (K
2^{2m})^{1-C}$. Consider an interval $[t_k, t_{k+1}]$. Clearly,
$|\tilde B_m(t_{k+1}) -  \tilde B_m(t_k)| = 2^{-m} \le 2^{-\half}
D^*_{K,m}$. On the other hand, $\tilde B_{m+1}$ makes 4 steps of
magnitude $2^{-m-1}$ on this interval. Then the maximum deviation
between $\tilde B_{m+1}$ and $\tilde B_m$ at the instant
$t_{k+\quart}$ (or at $t_{k+\frac34}$) cannot exceed $\frac{11}{4}
D^*_{K,m} + \quart 2^{-m} + 2^{-m-1} \le (\frac{11}{4} +
\frac{3}{4} 2^{-\half} ) D^*_{K,m}$. At time $t_{k+\half}$ the
deviation cannot be larger than $\frac{11}{4} D^*_{K,m} +
2^{-m-1}$, which is smaller than the previous bound. Hence
\[
\PP \lb \sup_{0\le t \le K} |\tilde{B}_{m+1}(t) -\tilde{B}_m(t)|
\ge \frac{11 + 3 \cdot 2^{-\half}}{4} D^*_{K,m} \rb \le 6
(K2^{2m})^{1-C} .
\]

Using the fact that
\[
\sum\limits_{j=0}^{\infty} (m+j)^{\frac34} 2^{-\frac{m+j}{2}} \le
m^{\frac34} 2^{-\frac{m}{2}} \sum\limits_{j=0}^{\infty}
(1+j)^{\frac34} 2^{-\frac{j}{2}} < \frac{65}{8} m^{\frac34}
2^{-\frac{m}{2}}
\]
(the last inequality can be checked for example by MAPLE), we
obtain that
\begin{eqnarray*}
\lefteqn {\PP \lb \sup_{j \ge 1} \sup_{0 \le t \le K} |B_{m+j}(t)
- B_m(t)| \ge \frac{65}{8} \cdot
\frac{11 + 3 \cdot 2^{-\half}}{4} D^*_{K,m} \rb } \\
& \le & \sum _{j=0}^{ \infty } \PP \lb \sup_{0 \le t \le K}
|B_{m+j+1}(t) - B_{m+j}(t)| \ge  \frac{11 + 3 \cdot 2^{-\half}}{4}
D^*_{K,m+j} \rb \\
& \le & \sum _{j=0}^{\infty } 6 (K 2^{2(m+j)})^{1-C} =
\frac{6}{1-4^{1-C}} (K 2^{2m})^{1-C} .
\end{eqnarray*}
Since $\frac{65}{8} \cdot \frac{11 + 3 \cdot 2^{-\half}}{4} < 27$,
this proves (b).

\end{proof}

It follows then that the above procedure gives a rather natural
and nearly optimal (as far as Skorohod embedding is concerned, see
below) construction of BM.

\begin{thm} \label{th:Wiener}
On any bounded interval the sequence $(\tilde B_m)$ almost surely
uniformly converges as $m \to \infty$ and the limit process is
Brownian motion $W$. For any $C>1$, and for any $K>0$ and $m \ge
1$ such that $K 2^{2m} \ge N_2(C)$, we have
\[
\PP \lb \sup_{0 \le t \le K} |W(t) - \tilde{B}_m(t)| \ge 27 \: C
K_*^{\quart} (\log_*K)^{\frac34} m^{\frac34} 2^{-\frac{m}{2}} \rb
\le \frac{6}{1-4^{1-C}} (K2^{2m})^{1-C} .
\]
\end{thm}

Now using the Borel--Cantelli lemma we get that for any fixed
$K>0$  there is a constant $c_K = 28 \; K_*^{\quart}
(\log_*K)^{\frac34}$ (taking $C = 1 + \frac{1}{27}$, say) such
that almost surely,
\begin{equation}\label{eq:Brownm}
\limsup_{m \to \infty} m^{-\frac34} 2^{\frac{m}{2}} \sup_{0 \le t
\le K} |W(t) - \tilde{B}_m(t)| <  c_K .
\end{equation}
Similarly, for any fixed $m \ge 1$ there is a constant $c_m = 55
\; m^{\frac34} 2^{-\frac{m}{2}}$ such that almost surely,
\begin{equation}\label{eq:BrownK}
\limsup_{K \to \infty} K^{-\quart} (\log K)^{-\frac34} \sup_{0 \le
t \le K} |W(t) - \tilde{B}_m(t)| < c_m .
\end{equation}
For, since $K^{\quart} (\log K)^{\frac34}$ and also for any
$\omega \in \Omega$, $\sup_{0 \le t \le K} |W(t) -
\tilde{B}_m(t)|$ are non-decreasing in $K$, it is enough to show
that
\begin{equation}\label{eq:asympt}
 \limsup_{n \to \infty} n^{-\quart} (\log n)^{-\frac34} \sup_{0 \le
t \le n+1} |W(t) - \tilde{B}_m(t)| < c_m ,
\end{equation}
when $n$ runs through integer values only. Since
\[
\frac{(n+1)^{\quart} (\log (n+1))^{\frac34}}{n^{\quart} (\log
n)^{\frac34}} \le 1 + \epsilon
\]
for any $\epsilon > 0$ if $n$ is large enough, taking, say,
$C=2+\frac{1}{27}$ in Theorem \ref{th:Wiener}, the Borel-Cantelli
lemma shows (\ref{eq:asympt}).

Next we are going to discuss the properties of another nested
sequence of random walks, obtained by Skorohod embedding. This
sequence is not identical, though asymptotically equivalent to the
above RW construction, cf. \cite[Theorem 4]{Sza96}. Given a Wiener
process $W$, first we define the stopping times which yield the
Skorohod embedded process $B_m(k2^{-2m})$ into $W$. For every
$m\ge 0$ let $s_m(0)=0$ and
\begin{equation}
s_m(k+1)=\inf{}\{s: s > s_m(k), |W(s)-W(s_m(k))|=2^{-m}\} \qquad
(k \ge 0). \label{eq:Skor1}
\end{equation}
With these stopping times the embedded process by definition is
\begin{equation}
B_m(k2^{-2m}) = W(s_m(k)) \qquad (m\ge 0, k\ge 0).
\label{eq:Skor2}
\end{equation}
This definition of $B_m$ can be extended to any real $t \ge 0$ by
pathwise linear interpolation. The next lemma describes some
useful facts about the relationship between $\tilde{B}_m$ and
$B_m$. These follow from \cite[Lemmas 5,7 and Theorem 4]{Sza96},
with some minor modifications.

In general, roughly saying, $\tilde{B}_m$ is more useful when
someone wants to generate stochastic processes from scratch, while
$B_m$ is more advantageous when someone needs a discrete
approximation of given processes, like in the case of stochastic
integration.

\begin{lem} \label{le:equid}
For any $C > 1$, and for any $K>0$ and $m \ge 1$ such that $K
2^{2m} \ge N_3(C)$ take the following subset of the sample space:
\[
 A^*_{K,m} = \lb \sup_{n > m} \: \sup_{0\le k2^{-2m}\le K}
|2^{-2n} T_{m,n}(k) - k2^{-2m}| < (42 \: C K \log_*K)^{\half}
m^{\half} 2^{-m} \rb ,
\]
where $T_{m,n}(k) = T_n \circ T_{n-1} \circ \cdots \circ T_m(k)$
for $n>m\ge 0$ and $k\ge 0$. Then
\[
\PP \lb (A^*_{K,m})^c \rb \le \frac{2}{1 - 4^{1-C}}(K2^{2m})^{1-C}
.
\]
Moreover, $\lim_{n\to \infty} 2^{-2n} T_{m,n}(k) = t_m(k)$ exists
almost surely and on $A^*_{K,m}$ we have
\[
\tilde{B}_m(k2^{-2m}) = W(t_m(k)) \qquad (0 \le k2^{-2m} \le K) ,
\]
cf. (\ref{eq:Skor2}). Further, on $A^*_{K,m}$ except for a zero
probability subset, $s_m(k) = t_m(k)$ and
\[
\sup_{0\le k2^{-2m}\le K} |s_m(k)-k2^{-2m}| \le (42 \: C K \log_*
K) ^{\half} m^{\half} 2^{-m} .
\]
\end{lem}

If the Wiener process is built by the RW construction described
above using a sequence $(\tilde{B}_m)$ of nested RW's and then one
constructs the Skorohod embedded RW's $(B_m)$, it is natural to
ask about rate of convergence of the latter. The answer described
by the next theorem is that it is essentially the same as the one
of $(\tilde{B}_m)$, cf. Theorem \ref{th:Wiener}.

\begin{thm} \label{th:Wienerm}
For any $C > 1$, and for any $K>0$ and $m \ge 1$ such that $K
2^{2m} \ge \max\{N_2(C),$ $N_3(C)\}$ we have
\[
\PP \lb \sup_{0\le t\le K} \left| W(t) - B_m(t) \right| \ge 27 \:
C K_*^{\quart} (\log_*K)^{\frac34} m^{\frac34} 2^{-\frac{m}{2}}
\rb \le \frac{8}{1 - 4^{1-C}} (K2^{2m})^{1-C} .
\]
\end{thm}
This theorem and its proof are slight modifications of Theorem 1
in \cite{SzeSza04}. Kiefer \cite{Kie69} showed that the best
possible rate of convergence one can get with Skorohod embedding
is $n^{\frac14} (\log n)^{\frac12} (\log \log n)^{\frac14}$. Our
rate of convergence $n^{\frac14} (\log n)^{\frac34}$, $n=K
2^{2m}$, is close to this.

\section{An elementary definition of Brownian local time}

We define the \emph{local time} of the random walk
$(\tilde{S}_m(k))_{k=0}^{\infty}$ at a point $x \in \mathbf{Z}$ at
time $k \in \mathbf{N}_0 = \{0,1,2,\dots\}$ as $\ell_m(0,x) = 0$
and
\[
\ell_m(k,x) = \#\{j : 0 \le j < k, \tilde{S}_m(j) = x \} \quad (k
\ge 1) .
\]
This is somewhat different from the more usual definition
$\tilde{\ell}_m(k,x) = \#\{j : 0 < j \le k, \tilde{S}_m(j) = x
\}$, but the former fits better the construction of Brownian
motion discussed in this paper. The local time of the $m$th
approximation $\tilde B_m$ at a point $x \in 2^{-m} \mathbf{Z}$ at
time $t \in 2^{-2m} \mathbf{N}_0$ is defined as $\loc_m(t,x) =
2^{-m} \ell_m \left(t2^{2m}, x 2^m \right)$, corresponding to the
fact that the spatial step size of $\tilde B_m$ is $2^{-m}$. This
is in complete agreement with (\ref{eq:Levyloctime}) replacing $W$
by $\tilde B_m$ there.

Finally, we define $\loc_m(t,x)$ for arbitrary $t \in
\mathbb{R}_+$ and $x \in \mathbb{R}$ by linear interpolation,
making it into a continuous process:
\begin{eqnarray*}
\loc_m(t,x) = \loc_m(t_k,x_j) &+& \frac{t-t_k}{t_{k+1}-t_k}
\left(\loc_m(t_{k+1},x_{j+1}) - \loc_m(t_k,x_{j+1}) \right) \\
&+& \frac{x-x_j}{x_{j+1}-x_j} \left(\loc_m(t_{k},x_{j+1}) -
\loc_m(t_k,x_{j}) \right)  \\
\end{eqnarray*}
if $x \ge x_j + 2^m (t-t_k)$,
\begin{eqnarray*}
\loc_m(t,x) = \loc_m(t_k,x_j) &+& \frac{t-t_k}{t_{k+1}-t_k}
\left(\loc_m(t_{k+1},x_{j}) - \loc_m(t_k,x_{j}) \right) \\
&+& \frac{x-x_j}{x_{j+1}-x_j} \left(\loc_m(t_{k+1},x_{j+1}) -
\loc_m(t_{k+1},x_{j}) \right)  \\
\end{eqnarray*}
if $x < x_j + 2^m (t-t_k)$, where $t_k = k 2^{-2m}$ $(k \in
\mathbf{N}_0)$, $x_j = j 2^{-m}$ $(j \in \mathbf{Z})$, and $(t,x)
\in [t_k, t_{k+1}] \times [x_j, x_{j+1}]$. Our aim is to define
the local time $\loc(t,x)$ of Brownian motion as the limit of
$\loc_m(t,x)$ as  $m \to \infty$ and to show that this limit is
jointly continuous in $(t,x)$.

The local times at $x=0$ are simply denoted by $\ell_m(k)$,
$\tilde{\ell}_m(k)$, and $\loc_m(t)$, respectively. Then
$\ell_m(k) = 1 + \tilde{\ell}_m(k-1)$ for $k \ge 1$ integer.

One can define ``one-sided'', \emph{up and down local times}
$\ell^+_m(k,x)$ and $\ell^-_m(k,x)$  ($k \in \mathbf{N}_0$, $x \in
\mathbf{Z}$) as well: $\ell^{\pm}_m(0,x) = 0$ and
\[
\ell^{\pm}_m(k,x) = \#\{j : 0 \le j < k, \: \tilde{S}_m(j) = x, \:
\tilde{S}_m(j+1) = x \pm 1 \}  \quad  (k \ge 1).
\]
Then $\ell_m(k,x) = \ell^+_m(k,x) + \ell^-_m(k,x)$. The
definitions of $\loc^+_m(t,x)$ and $\loc^-_m(t,x)$ for $t \in
\mathbb{R}_+$ and $x \in \mathbb{R}$ can go in the same way as the
definition of $\loc_m(t,x)$ above.

It can be useful to introduce an even finer division of local
time, introducing \emph{up-crossing local time} by
\[
\ell^{++}_m(k,x) = \#\{j : 1 \le j < k, \: \tilde{S}_{m}(j-1) =
x-1, \: \tilde{S}_{m}(j) = x, \: \tilde{S}_{m}(j+1) = x+1 \}
\]
and \emph{up-bouncing local time} by
\[
\ell^{-+}_{m}(k,x) = \#\{j : 1 \le j < k, \: \tilde{S}_{m}(j-1) =
x+1, \: \tilde{S}_{m}(j) = x, \: \tilde{S}_{m}(j+1) = x+1 \},
\]
where $k \ge 2$, $x \in \mathbb{Z}$. Then $\ell^+_{m}(k,x) =
\ell^{++}_{m}(k,x) + \ell^{-+}_{m}(k,x)$. The definitions of
\emph{down-crossing} and \emph{down-bouncing} local times
$\ell^{--}_m(k,x)$ and $\ell^{+-}_m(k,x)$ are similar.

As it is well-known, see e.g. \cite[p. 95]{Rev90}, the exact
distribution of the local time $\tilde{\ell}$ of a simple,
symmetric RW $(S(k))_{k=0}^{\infty}$, $S(0) = 0$, is
\[
\PP \lb \tilde{\ell}(2k) = j \rb = \PP \lb \tilde{\ell}(2k+1) = j
\rb = 2^{-2k+j} \binom{2k-j}{k} \quad (j=0,\dots,k) .
\]
Hence the usual argument for the De Moivre--Laplace limit theorem
gives that uniformly for any $0 \le j \le K_k = o(k^{\frac23})$
one has
\[
\PP \lb \tilde{\ell}(k) = j \rb \sim \sqrt{\frac{2}{\pi k}}
\exp\left( -\frac{j^2}{2k} \right) \qquad (k \to \infty),
\]
where $\sim$ denotes ``asymptotically equal''. Then for any
sequence $u_k \to \infty$, $u_k = o(k^{\frac16})$ we obtain the
following large deviation inequality for local times $\ell(k,x)$,
$\ell^+(k,x)$, and $\ell^-(k,x)$ $(x \in \mathbb{Z})$:
\begin{equation}\label{eq:loctasy}
\PP \lb \frac{\ell^{\pm}(k,x)}{\sqrt{k}} \ge u_k \rb \le \PP \lb
\frac{\ell(k,x)}{\sqrt{k}} \ge u_k \rb \le \PP \lb
\frac{\ell(k)}{\sqrt{k}} \ge u_k \rb \le
\exp\left(-\frac{u_k^2}{2} \right)
\end{equation}
if $k$ is large enough, $k \ge k_0$.

Both the statements and the proofs of the next lemmas about local
times are very similar to the ones about BM approximations $\tilde
B_m(t)$ in Lemma \ref{le:refin}.

\begin{lem} \label{le:loctime1}
For any $C>1$, and for any $K>0$ and $m \ge 1$ such that $K 2^{2m}
\ge N_4(C)$, we have
\begin{eqnarray*}
\PP \lb \sup_{j \in \mathbb{Z}} \sup_{0\le t_k \le K}
|\loc_{m+1}(t_k, x_j) -\loc_m(t_k, x_j)| \ge 6 \: C K^{\quart}
(\log_*K)^{\frac34} m^{\frac34}  2^{-\frac{m}{2}} \rb \\
\le 12 (K2^{2m})^{1-C} ,
\end{eqnarray*}
where $t_k = k 2^{-2m}$ and $x_j = j 2^{-m}$.
\end{lem}

\begin{proof}
We are making use of the fact that a RW of length $N = K 2^{2m}$
typically cannot have values or local times much larger than
$\sqrt{N}$. Let us introduce the abbreviations $D_{K,m} = C
K^{\quart} (\log_*K)^{\frac34} m^{\frac34} 2^{-\frac{m}{2}}$, $M =
(3 C m K \log_*K)^{\half} 2^m$ (for the ``maximum''), $N'' =
\sqrt{3} \: M$ (for the ``maximal'' local time), and $K_m =
\lfloor K 2^{2m} \rfloor$. Then using the triangle inequality
\begin{eqnarray*}
\lefteqn{|\ell_{m+1}(4k,2j) - 2 \ell_{m}(k,j)|} \\
&\le& |\ell_{m+1}(T_{m+1}(k),2j) - 2 \ell_{m}(k,j)| +
|\ell_{m+1}(4k,2j) - \ell_{m+1}(T_{m+1}(k),2j)| \\
&=:& A_{m+1}(k,j) + B_{m+1}(k,j) ,
\end{eqnarray*}
we get that
\begin{eqnarray*}
\lefteqn{\PP \lb \sup_{j \in \mathbb{Z}} \sup_{0 \le t_k \le K}
|\loc_{m+1}(t_k, x_j) - \loc_m(t_k, x_j)| \ge 6 D_{K,m} \rb} \\
&\le& \PP \lb \sup_{|j| \le M} \sup_{0 \le t_k \le K} 2^{-m-1}
|\ell_{m+1}(4k,2j) - 2 \ell_{m}(k,j)| \ge 6 D_{K,m} \rb \\
&+& \PP \lb \sup_{|j| > M} \sup_{0 \le t_k \le K}
|\loc_{m+1}(t_k, x_j) - \loc_m(t_k, x_j)| \ge 6 D_{K,m} \rb \\
&\le& \sum\limits_{|j| \le M} \sum\limits_{n=1}^{N''} \PP \lb
\ell_m(K_m,j) = n; \sup_{0 \le k \le K_m} 2^{-m-1}
 A_{m+1}(k,j) \ge 3 D_{K,m} \rb \\
&+& \sum\limits_{|j| \le M} \sum\limits_{n=N''+1}^{K_m /2} \PP \lb
\ell_m(K_m,j) = n; \sup_{0 \le k \le K_m}  2^{-m-1}
A_{m+1}(k,j) \ge 3 D_{K,m} \rb \\
&+& \PP \lb \sup_{|j| \le M} \sup_{0 \le k \le K_m} 2^{-m-1}
B_{m+1}(k,j)
\ge 3 D_{K,m} \rb \\
&+&  \PP \lb \sup_{|j| > M} \sup_{0 \le t_k \le K}
|\loc_{m+1}(t_k, x_j) - \loc_m(t_k, x_j)| \ge 6 D_{K,m} \rb \\
&=:& p_1 + p_2 + p_3 + p_4 ,
\end{eqnarray*}
We are going to estimate $p_1$, $p_2$, $p_3$, and $p_4$ one by
one.

\emph{Estimation of $p_1$.}

This is the essential part of the proof. Fixing $m \ge 0$, $j \in
[-M, M]$, $n \ge 1$, and given $\ell_m(K_m,j) = n$, let $\tau_0  <
\tau_1 < \cdots < \tau_{n-1} < K_m$ denote the random time
instants when the ``twisted'' RW $\tilde S_m(k)$ is at the point
$j$ in the interval $[0, K_m)$: $\tilde S_m(\tau_i) = j$; then
$\tilde S_{m+1}(T_{m+1}(\tau_i)) = 2j$. We define $\tau_n = K_m$,
irrespective whether $\tilde S_m(\tau_n) = 0$ or not. Then we can
write that
\begin{eqnarray*}
\ell_{m+1}(T_{m+1}(K_m),2j) &=& \sum\limits_{i=1}^{n}
\ell_{m+1}(T_{m+1}(\tau_{i}),2j)  -
\ell_{m+1}(T_{m+1}(\tau_{i-1}),2j) \\
&=& \sum\limits_{i=1}^{n} \gamma_{n,i} .
\end{eqnarray*}
(The dependence of $\gamma_{n,i}$ on $m$ and $j$ is suppressed in
the notation.) Here each random variable $\gamma_{n,i}$ ($i = 1,
\dots , n$) is the number of time instants when $\tilde
S_{m+1}(k)$ hits the point $2j$ in the interval
$\left[T_{m+1}(\tau_{i-1}), T_{m+1}(\tau_{i})\right)$. This is
simply $1$ plus the number of the $+1, -1$ or $-1, +1$ pairs of
steps of $\tilde S_{m+1}(k)$, starting from $T_{m+1}(\tau_{i-1})$.
Any such sequence ends with a pair $+1, +1$ or $-1, -1$. Clearly,
this means that $(\gamma_{n,i})_{i=1}^{n}$ is a sequence of
independent, geometrically distributed random variables with
parameter $p = \half$; $\EE (\gamma_{n,i}) = 2$ and $\Var
(\gamma_{n,i}) = 2$.

From the construction of $\tilde S_{m+1}(k)$ discussed in Section
2 it is also clear that given $\ell_m(K_m,j) = n$ $(n \ge 1)$,
\begin{eqnarray*}
\lefteqn{\sup_{0 \le k \le K_m} |\ell_{m+1}(T_{m+1}(k),2j) - 2
\ell_{m}(k,j)|} \\
&=& \sup_{1 \le i \le n} |\ell_{m+1}(T_{m+1}(\tau_i),2j) - 2
\ell_{m}(\tau_i,j)| = \sup_{1 \le r \le n} \left|
\sum\limits_{i=1}^{r} \gamma_{n,i} - 2r \right| .
\end{eqnarray*}
(Note that $\ell_m(K_m,j) = 0$ implies
$\ell_{m+1}(T_{m+1}(K_m),2j) = 0$ as well.) Hence
\begin{eqnarray*}
p_1 &\le& \sum\limits_{|j| \le M} \sum\limits_{n=1}^{N''} \PP \lb
\sup_{1 \le r \le n} \left| \sum\limits_{i=1}^{r} \gamma_{n,i} -
2r \right| \ge 6 \: 2^m D_{K,m} \rb \\
&\le& (2 M + 1) N'' \: \PP \lb \sup_{1 \le r \le N''} \left|
\sum\limits_{i=1}^{r} \frac{\gamma_{i} - 2}{\sqrt{2}} \right| \ge
3\sqrt{2} \: 2^m D_{K,m} \rb
\end{eqnarray*}

Here we are going to use Lemma \ref{le:ldi} with $N''$ and $C''$.
The constant $C''$ is chosen so that the error probabilities
$p_1$, $p_2$, $p_3$, and $p_4$ be of the same order: $(2 M + 1)
N'' (N'')^{1-C''} \le (2.5/\sqrt{3}) (N'')^{3-C''} \approx
N^{1-C}$, where $N=K 2^{2m}$. Then $(N'')^{3-C''} \le (K^{\half}
2^m)^{3-C''} = (K 2^{2m})^{1-C}$ implies $C'' = 2C+1 < 3C$. Since
$\log N'' \le m \log_*C \log_*K$ if $K$ or $m$ are large enough,
it follows that $(2 C'' N'' \log N'')^{\half} \le 3\sqrt{2}
\:2^{m} D_{K,m}$. Thus we obtain that
\begin{equation}\label{eq:p1}
p_1 \le (5 / \sqrt{3}) (N'')^{3-C''} < 3 (K 2^{2m})^{1-C} ,
\end{equation}
if $K 2^{2m}$ is large enough.

\emph{Estimation of $p_2$.}

Here we are using inequality (\ref{eq:loctasy}):
\begin{eqnarray}
p_2 &\le& \sum\limits_{|j| \le M} \sum\limits_{n=N''+1}^{K_m/2}
\PP \lb \ell_m(K_m, j) = n \rb \le (2 M + 1) \PP \lb \ell_m(K_m,
0) \ge N''
\rb \nonumber \\
&\le& (2 M + 1) \PP \lb \frac{\ell_m(K_m)}{\sqrt{K_m}} \ge (9 C m
\log_*K)^\half \rb \nonumber \\
&\le& (2 M + 1) \exp\left(-\frac92 C m \log_*K\right)
\nonumber \\
&\le& (K 2^{2m})^{1-C} , \label{eq:p2}
\end{eqnarray}
if $K 2^{2m}$ is large enough.

\emph{Estimation of $p_3$.}

Here Lemma \ref{le:tlag} and inequality (\ref{eq:loctasy}) will be
used. Let $A_{K,m}$ and $N'$ be the same as in (\ref{eq:AKm}).
Then
\begin{eqnarray*}
p_3 &\le& \PP \lb A_{K,m} \rb \nonumber \\
&+& \sum\limits_{|j| \le M} \sum\limits_{k=0}^{K 2^{2m}} \PP \lb
A^c_{K,m} ; |\ell_{m+1}(4k,2j) - \ell_{m+1}(T_{m+1}(k),2j)|
\ge 6 \: 2^m  D_{K,m} \rb \\
&\le& 2 (K 2^{2m})^{1-C} \nonumber \\
&+& \sum\limits_{|j| \le M} \sum\limits_{k=0}^{K 2^{2m}} \PP \lb
\sup_{|i-4k| \le N'} |\ell_{m+1}(4k,2j) - \ell_{m+1}(i,2j)| \ge 6
\: 2^m D_{K,m} \rb .
\end{eqnarray*}
If $i \ge 4k$, say, then
\[
\PP \lb \ell_{m+1}(i,2j) - \ell_{m+1}(4k,2j) \ge u \rb \le \PP \lb
\ell_{m+1}(i-4k,0) \ge u \rb ,
\]
for any $u \ge 0$. For, $\tilde S_{m+1}(0) = 0$, while $\tilde
S_{m+1}(4k)$ can be different from $2j$ and after the first visit
$\tau \ge 4k$ of $\tilde S_{m+1}$ to the point $2j$ it starts from
scratch by the strong Markov property of a simple, symmetric RW.
We get that
\begin{eqnarray}
p_3 &\le& 2 (K 2^{2m})^{1-C} + 2 (2M+1) K 2^{2m} \PP \lb \sup_{0
\le i-4k \le N'} \ell_{m+1}(i-4k) \ge 6 \: 2^m D_{K,m} \rb
\nonumber \\
&\le& 2 (K 2^{2m})^{1-C} + 2 (2M+1) K 2^{2m} \PP \lb
\frac{\ell_{m+1}(N')}
{\sqrt{N'}} \ge 2^\quart (3C)^{\frac34} (m \log_*K)^\half \rb
\nonumber \\
&\le& 2 (K 2^{2m})^{1-C} + 2 (2M+1) K 2^{2m} \exp
\left(-2^{-\half} (3C)^{\frac32} m \log_*K \right) \nonumber \\
&\le& 4 (K 2^{2m})^{1-C} , \label{eq:p3}
\end{eqnarray}
if $K 2^{2m}$ is large enough.

\emph{Estimation of $p_4$.}

Here we are going to apply inequality (\ref{eq:ldp2}) to
$\tilde{S}_{m+1}(4k)$ and $\tilde{S}_{m}(k)$:
\begin{eqnarray}
p_4 &\le& \PP \lb \sup_{0 \le t_k \le K} |\tilde{B}_{m+1}(t_k)|
\ge M 2^{-m} \rb + \PP \lb \sup_{0 \le t_k \le K}
|\tilde{B}_{m}(t_k)|
\ge M 2^{-m} \rb \nonumber \\
&=& \PP \lb \sup_{0 \le k \le K_m} |\tilde{S}_{m+1}(4k)| \ge 2 M
\rb + \PP \lb \sup_{0 \le k \le K_m} |\tilde{S}_{m}(k)|
\ge M  \rb \nonumber \\
&\ge& \PP \lb \sup_{0 \le k \le K_m} |\tilde{S}_{m+1}(4k)| \ge
\left(2 C K 2^{2(m+1)} \log(K 2^{2m})\right)^{\half} \rb
\nonumber \\
&& + \: \PP \lb \sup_{0 \le k \le K_m} |\tilde{S}_{m}(k)| \ge
\left(2 C K 2^{2m} \log(K 2^{2m})\right)^{\half} \rb \nonumber \\
&\le& 4 (K 2^{2m})^{1-C} , \label{eq:p4}
\end{eqnarray}
if $K 2^{2m}$ is large enough.

Combining (\ref{eq:p1}), (\ref{eq:p2}), (\ref{eq:p3}), and
(\ref{eq:p4}), we get the statement of the lemma.
\end{proof}

\begin{lem} \label{le:loctime2}
For any $C>1$, and for any $K>0$ and $m \ge 1$ such that $K 2^{2m}
\ge N_5(C)$, we have
\begin{eqnarray*}
\PP \lb \sup_{j \in \mathbb{Z}} \sup_{0\le t_k \le K}
|\loc^{+}_{m+1}(t_k, x_j) - \half \loc_m(t_k, x_j)| \ge 6 \: C
K^{\quart} (\log_*K)^{\frac34} m^{\frac34}  2^{-\frac{m}{2}} \rb \\
\le 12 (K2^{2m})^{1-C} ,
\end{eqnarray*}
where $t_k = k 2^{-2m}$ and $x_j = j 2^{-m}$. A similar statement
holds for $\loc^{-}_{m+1}(t_k, x_j)$.
\end{lem}

\begin{proof}
This proof goes similarly as the proof of the previous Lemma
\ref{le:loctime1}, except for the estimation of $p_1$, where there
are some differences. Hence this is the only point detailed in the
sequel. We only discuss the case of $\loc^+_{m+1}$ as the case of
$\loc^-_{m+1}$ is analogous.

Fixing $m \ge 0$, $j \in [-M, M]$, $n \ge 1$, and given
$\ell_m(K_m,j) = n$, let $\tau_0  < \tau_1 < \cdots < \tau_{n-1} <
K_m$ denote the random time instants when $\tilde S_m(\tau_i) = j$
in the interval $[0, K_m)$; then $\tilde S_{m+1}(T_{m+1}(\tau_i))
= 2j$. We define $\tau_n = K_m$. Then we can write that
\begin{eqnarray*}
\ell^+_{m+1}(T_{m+1}(K_m),2j) &=& \sum\limits_{i=1}^{n}
\ell^+_{m+1}(T_{m+1}(\tau_{i}),2j)  -
\ell^+_{m+1}(T_{m+1}(\tau_{i-1}),2j) \\
&=& \sum\limits_{i=1}^{n} (\alpha_{n,i} + X_{n,i}) .
\end{eqnarray*}
(The dependence of $\alpha_{n,i}$ and $X_{n,i}$ on $m$ and $j$ is
suppressed in the notation.) Here each random variable
$\alpha_{n,i}$ ($i = 1, \dots , n$) is the number of time instants
$k \in \left[T_{m+1}(\tau_{i-1}), T_{m+1}(\tau_{i})\right)$ when
$\tilde S_{m+1}(k) = 2j$, $\tilde S_{m+1}(k+1) = 2j+1$, and
$\tilde S_{m+1}(k+2) = 2j$. This is simply the number of the $+1,
-1$ pairs of steps of $\tilde S_{m+1}(k)$ in the interval
$\left[T_{m+1}(\tau_{i-1}), T_{m+1}(\tau_{i})\right)$. Any
sequence of $+1, -1$ or $-1, +1$ pairs ends with a pair $+1, +1$
or $-1, -1$. In the former case, $X_{n,i} = 1$, in the second case
$X_{n,i} = 0$. It follows that $(\alpha_{n,i})_{i=1}^{n}$ is a
sequence of independent, geometrically distributed random
variables with parameter $p = \frac23$; $\EE (\alpha_{n,i}) =
\half$ and $\Var (\alpha_{n,i}) = \frac34$. Further,
$(X_{n,i})_{i=1}^{n}$ is a sequence of independent indicator
variables with parameter $p = \half$; $\EE (X_{n,i}) = \half$ and
$\Var (X_{n,i}) = \frac14$. The two sequences are also
independent.

It is also clear that given $\ell_m(K_m,j) = n$ $(n \ge 1)$,
\begin{eqnarray*}
\lefteqn{C_{m+1}(k,j) := \sup_{0 \le k \le K_m}
|\ell^+_{m+1}(T_{m+1}(k),2j) -
\ell_{m}(k,j)|} \\
&=& \sup_{1 \le i \le n} |\ell^+_{m+1}(T_{m+1}(\tau_i),2j) -
\ell_{m}(\tau_i,j)| = \sup_{1 \le r \le n} \left|
\sum\limits_{i=1}^{r} (\alpha_{n,i} + X_{n,i}) - r \right| .
\end{eqnarray*}
Hence in the same way as in the proof of Lemma \ref{le:loctime1},
\begin{eqnarray}
p_1 &:=& \sum\limits_{|j| \le M} \sum\limits_{n=1}^{N''} \PP \lb
\ell_m(K_m,j) = n;  2^{-m-1} C_{m+1}(k,j) \ge 3 D_{K,m} \rb
\nonumber \\
&\le& \sum\limits_{|j| \le M} \sum\limits_{n=1}^{N''} \PP \lb
\sup_{1 \le r \le n} \left| \sum\limits_{i=1}^{r} (\alpha_{n,i} +
X_{n,i}) - r \right| \ge 6 \: 2^m D_{K,m} \rb \nonumber \\
&\le& (2 M + 1) N'' \: \PP \lb \sup_{1 \le r \le N''} \left|
\sum\limits_{i=1}^{r} (\alpha_{i} + X_{n,i} - 1) \right|
\ge 6 \: 2^m D_{K,m} \rb \nonumber \\
&<& 3 (K 2^{2m})^{1-C}  \label{eq:p1+},
\end{eqnarray}
if $K 2^{2m}$ is large enough. This ends the proof of the lemma.

\end{proof}

\begin{lem} \label{le:loctime3}
For any $C>1$, and for any $K>0$ and $m \ge 1$ such that $K 2^{2m}
\ge N_6(C)$, we have
\begin{eqnarray*}
\PP \lb \sup_{r \ge 1} \sup_{(t,x) \in [0,K] \times \mathbb{R}}
|\loc_{m+r}(t,x) -\loc_m(t,x)| \ge 79 \: C K_*^{\quart}
(\log_*K)^{\frac34}
m^{\frac34}  2^{-\frac{m}{2}} \rb \\
\le \frac{15}{1-4^{1-C}} (K2^{2m})^{1-C} ,
\end{eqnarray*}
where $K_* = \max\{1, K\}$.
\end{lem}

\begin{proof}

{\sc Step 1}

$\tilde{B}_{m+1}(k 2^{-2(m+1)})$ visits new points $x_{j+\half} =
(j + \half) 2^{-m} = (2j+1) 2^{-m-1}$ $(j \in \mathbb{Z})$ that
were not visited by $\tilde{B}_{m}(k 2^{-2m})$. Thus first we need
to show an inequality for $\loc_{m+1}(t_{k},x_{j+\half})$, similar
to the one in Lemma \ref{le:loctime1}:
\begin{eqnarray}
\PP \lb \sup_{j \in \mathbb{Z}} \sup_{0 \le t_k \le K}
|\loc_{m+1}(t_k, x_{j+\half}) - \loc_m(t_k, x_{j+\half})| \ge 9 \:
C K^{\quart} (\log_*K)^{\frac34} m^{\frac34}  2^{-\frac{m}{2}} \rb
\nonumber \\
\le 15 (K2^{2m})^{1-C} . \label{eq:jhalf}
\end{eqnarray}
Since $\loc_m(t_k, x_{j+\half})$ is obtained by linear
interpolation and
\[
\ell_{m+1}(T_{m+1}(k),2j+1) = \ell^+_{m+1}(T_{m+1}(k),2j) +
\ell^-_{m+1}(T_{m+1}(k),2j+2)
\]
(note that $\tilde{S}_{m+1}(T_{m+1}(k))$ never equals $2j+1$), it
follows that
\begin{eqnarray*}
\lefteqn{|\loc_{m+1}(t_k, x_{j+\half}) - \loc_m(t_k, x_{j+\half})|}
\\
&=& 2^{-m-1} |\ell_{m+1}(4k,2j+1) - \ell_{m}(k,j) - \ell_{m}(k,j+1)|
\\
&\le& 2^{-m-1} |\ell_{m+1}(T_{m+1}(k),2j+1) - \ell_{m}(k,j) -
\ell_{m}(k,j+1)| \\
&+& 2^{-m-1} |\ell_{m+1}(4k,2j+1) - \ell_{m+1}(T_{m+1}(k),2j+1)|
\\
&\le& 2^{-m-1} |\ell^+_{m+1}(T_{m+1}(k),2j) - \ell_{m}(k,j)| \\
&+& 2^{-m-1} |\ell^-_{m+1}(T_{m+1}(k),2j+2) - \ell_{m}(k,j+1)| \\
&+& 2^{-m-1} |\ell_{m+1}(4k,2j+1) - \ell_{m+1}(T_{m+1}(k),2j+1)| \\
&=:& F_{m+1}(k,j) + G_{m+1}(k,j) + H_{m+1}(k,j).
\end{eqnarray*}

From this point the proof of (\ref{eq:jhalf}) goes similarly as
the proof of Lemma \ref{le:loctime1}, except for the estimation of
$p_1$, where there are some differences. Hence this is the only
point detailed here. By (\ref{eq:p1+}),
\begin{eqnarray*}
p_1 &:=& \sum\limits_{|j| \le M} \sum\limits_{n=1}^{N''} \PP \lb
\ell_m(K_m,j) = n; \sup_{0 \le k \le K_m} F_{m+1}(k,j) \ge 3
D_{K,m} \rb
\\
&+& \sum\limits_{|j| \le M} \sum\limits_{n=1}^{N''} \PP \lb
\ell_m(K_m,j+1) = n; \sup_{0 \le k \le K_m} G_{m+1}(k,j) \ge 3
D_{K,m}
\rb \\
&<& 6 (K 2^{2m})^{1-C},
\end{eqnarray*}
if $K 2^{2m}$ is large enough.

{\sc Step 2}

Let $D^*_{K,m} = C K_*^{\quart} (\log_*K)^{\frac34} m^{\frac34}
2^{-\frac{m}{2}}$. By Lemma  \ref{le:loctime1},
\[
\sup_{j \in \mathbb{Z}} \sup_{0\le t_k \le K} |\loc_{m+1}(t_k,x_j)
-  \loc_m(t_k,x_j)| < 6 D^*_{K,m} ,
\]
except for an event of probability not exceeding $12 (K
2^{2m})^{1-C}$. Also, by (\ref{eq:jhalf}),
\[
\sup_{j \in \mathbb{Z}} \sup_{0\le t_k \le K}
|\loc_{m+1}(t_k,x_{j+\half}) -  \loc_m(t_k,x_{j+\half})| < 9
D^*_{K,m} ,
\]
except for an event of probability not exceeding $15 (K
2^{2m})^{1-C}$.

Consider an interval $[t_k, t_{k+1}]$. Clearly, $0 \le
\loc_m(t_{k+1},x_j) - \loc_m(t_k,x_j) \le 2^{-m} \le 2^{-\half}
D^*_{K,m}$. On the other hand, $\tilde B_{m+1}$ makes 4 steps on
this interval. Thus $0 \le \loc_{m+1}(t_{k+1},x_j) -
\loc_{m+1}(t_k,x_j) \le 2 \cdot 2^{-m-1} = 2^{-m}$ and $0 \le
\loc_{m+1}(t_{k+1},x_{j+\half}) - \loc_{m+1}(t_k,x_{j+\half}) \le
2^{-m}$as well. Since $\loc_{m}(t,x)$ and $\loc_{m+1}(t,x)$ are
obtained by linear interpolation for real $t$ and $x$, it follows
that
\begin{equation}\label{eq:loctime3}
\PP \lb \sup_{(t,x) \in [0,K] \times \mathbb{R}} |\loc_{m+1}(t,x)
-\loc_m(t,x)| \ge (9 + 2^{-\half}) D^*_{K,m} \rb \le 15
(K2^{2m})^{1-C} .
\end{equation}

From this point the proof is the same as the last part of the
proof of Lemma \ref{le:refin}(b), except for the constant
multipliers. Since $\frac{65}{8} (9 + 2^{-\half}) < 79$, this
proves the statement of the lemma.

\end{proof}

By the Borel--Cantelli lemma, the previous Lemma \ref{le:loctime3}
leads to an alternative definition of the local time of Brownian
motion $W$ obtained in Theorem \ref{th:Wiener}, via a sequence of
continuous local times of the approximations $\tilde B_m$.

\begin{thm} \label{th:loctime}
On any strip $[0, K] \times \mathbb{R}$ the sequence
$(\loc_m(t,x))$ almost surely uniformly converges as $m \to
\infty$ and the limit is a process $\loc(t,x)$ jointly continuous
in $(t,x)$, the local time of Brownian motion $W(t)$. For any
$C>1$, and for any $K>0$ and $m \ge 1$ such that $K 2^{2m} \ge
N_6(C)$, we have
\begin{eqnarray*}
&\PP& \lb \sup_{(t,x) \in [0,K] \times \mathbb{R}} |\loc(t,x) -
\loc_m(t,x)| \ge 79 \: C K_*^{\quart}
(\log_*K)^{\frac34} m^{\frac34} 2^{-\frac{m}{2}} \rb \\
&&\le \frac{15}{1-4^{1-C}} (K2^{2m})^{1-C} .
\end{eqnarray*}
\end{thm}

Now using the Borel--Cantelli lemma we get that for any fixed
$K>0$  there is a constant $c_K = 80 \; K_*^{\quart}
(\log_*K)^{\frac34}$ (taking $C = 1 + \frac{1}{79}$, say) such
that almost surely,
\begin{equation}\label{eq:loctm}
\limsup_{m \to \infty} m^{-\frac34} 2^{\frac{m}{2}} \sup_{(t,x)
\in [0,K] \times \mathbb{R}} |\loc(t,x) - \loc_m(t,x)| <  c_K .
\end{equation}
Also, for any fixed $m \ge 1$ there is a constant $c_m = 159 \;
m^{\frac34} 2^{-\frac{m}{2}}$ (taking $C = 2 + \frac{1}{79}$, say)
such that almost surely,
\begin{equation}\label{eq:loctK}
\limsup_{K \to \infty} K^{-\quart} (\log K)^{-\frac34} \sup_{(t,x)
\in [0,K] \times \mathbb{R}} |\loc(t,x) - \loc_m(t,x)| < c_m .
\end{equation}
This also follows by the Borel-Cantelli lemma, in a similar way as
(\ref{eq:BrownK}) did.

One has similar convergence results for the one-sided local times
as well.

\begin{thm} \label{th:loctime+}
On any strip $[0, K] \times \mathbb{R}$ the sequence
$(\loc^+_m(t,x))$ almost surely uniformly converges as $m \to
\infty$ to the one half of the Brownian local time $\loc(t,x)$.
For any $C>1$, and for any $K>0$ and $m \ge 1$ such that $K 2^{2m}
\ge N_7(C)$, we have
\begin{eqnarray*}
&\PP& \lb \sup_{(t,x) \in [0,K] \times \mathbb{R}} \left|\half \:
\loc(t,x) - \loc^+_{m+1}(t,x)\right| \ge 50 \: C K_*^{\quart}
(\log_*K)^{\frac34} m^{\frac34} 2^{-\frac{m}{2}} \rb \\
&&\le \frac{30}{1-4^{1-C}} (K2^{2m})^{1-C} .
\end{eqnarray*}
Similar statements hold for $(\loc^-_{m+1}(t,x))$ as well.
\end{thm}
\begin{proof}

{\sc Step 1}

First we need to show an inequality similar to (\ref{eq:jhalf})
for points $x_{j + \half}$:
\begin{eqnarray}
&\PP& \lb\sup_{j \in \mathbb{Z}} \sup_{0 \le t_k \le K}
\left|\loc^+_{m+1}(t_k, x_{j+\half}) - \half \: \loc_m(t_k,
x_{j+\half})\right| \ge 9 \: C K^{\quart} (\log_*K)^{\frac34}
m^{\frac34} 2^{-\frac{m}{2}} \rb
\nonumber \\
&&\le 15 (K2^{2m})^{1-C} . \label{eq:jhalf2}
\end{eqnarray}
Our argument will follow a similar path to the ones in Lemmas
\ref{le:loctime2} and \ref{le:loctime3}. Since $\loc_m(t_k,
x_{j+\half})$ is obtained by linear interpolation, it follows that
\begin{eqnarray*}
\lefteqn{\left|\loc^+_{m+1}(t_k, x_{j+\half}) - \half
\:\loc_m(t_k, x_{j+\half})\right|}
\\
&=& 2^{-m-1} \left|\ell^+_{m+1}(4k,2j+1) - \half \left(
\ell_{m}(k,j) + \ell_{m}(k,j+1)\right)\right|
\\
&\le& 2^{-m-1} \left|\ell^+_{m+1}(T_{m+1}(k),2j+1) - \half
\left(\ell_{m}(k,j) +
\ell_{m}(k,j+1)\right) \right| \\
&+& 2^{-m-1} |\ell^+_{m+1}(4k,2j+1) -
\ell^+_{m+1}(T_{m+1}(k),2j+1)|
\\
&\le& 2^{-m-1} \left|\ell^{++}_{m+1}(T_{m+1}(k),2j+1) - \half \:
\ell_{m}(k,j)\right| \\
&+& 2^{-m-1} \left|\ell^{-+}_{m+1}(T_{m+1}(k),2j+1)
- \half \: \ell_{m}(k,j+1)\right| \\
&+& 2^{-m-1} |\ell^+_{m+1}(4k,2j+1) -
\ell^+_{m+1}(T_{m+1}(k),2j+1)| ,
\end{eqnarray*}
where we applied the notations for up-crossing and up-bouncing
local times, introduced above.

Here, analogously to the proof of Lemma \ref{le:loctime2}, given
$\ell_m(K_m,j) = n$, one has
\[
\ell^{++}_{m+1}(T_{m+1}(K_m),2j+1) =  \sum\limits_{i=1}^{n}
X_{n,i} ,
\]
where $X_{n,i} = 1$ if a sequence of $+1, -1$ or $-1, +1$ pairs of
steps of $\tilde S_{m+1}(k)$ in the interval
$\left[T_{m+1}(\tau_{i-1}), T_{m+1}(\tau_{i})\right)$ ends with a
pair $+1, +1$ and $0$ otherwise; $\tau_0  < \tau_1 < \cdots <
\tau_{n-1} < K_m$ are the random time instants when $\tilde
S_m(\tau_i) = j$ in the interval $[0, K_m)$ and $\tau_n = K_m$.
Then $(X_{n,i})_{i=1}^{n}$ is a sequence of independent indicator
variables with parameter $p = \half$; $\EE (X_{n,i}) = \half$ and
$\Var (X_{n,i}) = \frac14$. Further, given $\ell_m(K_m,j) = n$,
\[
\sup_{0 \le k \le K_m} \left|\ell^{++}_{m+1}(T_{m+1}(k),2j+1) -
\half \: \ell_{m}(k,j)\right| = \sup_{1 \le r \le n} \left|
\sum\limits_{i=1}^{r} \left(X_{n,i} - \half \right) \right| .
\]

Similarly, given $\ell_m(K_m,j+1) = n$,
\[
\ell^{-+}_{m+1}(T_{m+1}(K_m),2j+1) =  \sum\limits_{i=1}^{n}
\beta_{n,i} ,
\]
where $\beta_{n,i}$ is simply the number of the $-1, +1$ pairs of
steps of $\tilde S_{m+1}(k)$ in the interval
$\left[T_{m+1}(\tau_{i-1}), T_{m+1}(\tau_{i})\right)$; $\tau_0  <
\tau_1 < \cdots < \tau_{n-1} < K_m$ are the random time instants
when $\tilde S_m(\tau_i) = j+1$ in the interval $[0, K_m)$ and
$\tau_n = K_m$. It follows that $(\beta_{n,i})_{i=1}^{n}$ is a
sequence of independent, geometrically distributed random
variables with parameter $p = \frac23$; $\EE (\beta_{n,i}) =
\half$ and $\Var (\beta_{n,i}) = \frac34$. Moreover, given
$\ell_m(K_m,j+1) = n$,
\[
\sup_{0 \le k \le K_m} \left|\ell^{-+}_{m+1}(T_{m+1}(k),2j+1) -
\half \: \ell_{m}(k,j+1)\right| = \sup_{1 \le r \le n} \left|
\sum\limits_{i=1}^{r} \left(\beta_{n,i} - \half \right) \right| .
\]

From this point the proof of (\ref{eq:jhalf2}) is essentially the
same as the proof of (\ref{eq:jhalf}).

{\sc Step 2}

By Lemma \ref{le:loctime2} and formula (\ref{eq:jhalf2}), using
the same argument as in {\sc Step 2} of Lemma \ref{le:loctime3},
it follows that
\begin{equation}\label{eq:loctime+}
\PP \lb \sup_{(t,x) \in [0,K] \times \mathbb{R}}
\left|\loc^+_{m+1}(t,x) - \half \: \loc_m(t,x)\right| \ge (9 +
2^{-\half}) D^*_{K,m} \rb \le 15 (K2^{2m})^{1-C} .
\end{equation}
Since $|\half \: \loc(t,x) - \loc^+_{m+1}(t,x)| \le |\half \:
\loc(t,x) - \half \: \loc_m(t,x)| + |\half \: \loc_m(t,x) -
\loc^+_{m+1}(t,x)|$, formula (\ref{eq:loctime+}) and Theorem
\ref{th:loctime} give that
\begin{eqnarray*}
&\PP& \lb \sup_{(t,x) \in [0,K] \times \mathbb{R}} \left|\half \:
\loc(t,x) - \loc^+_{m+1}(t,x) \right| \ge \left(\frac{79}{2} + 9 +
2^{-\half}\right) D^*_{K,m} \rb \\
&& \le \left(\frac{15}{1-4^{1-C}} + 15\right) (K2^{2m})^{1-C} .
\end{eqnarray*}
This proves the theorem.
\end{proof}

We mention that similar convergence results can be shown for up-
and down-crossing or up- and down-bouncing local times as well.


\end{document}